\documentclass[12pt,draft]{extartic}

\usepackage{amsmath,amsthm,amssymb,calc}
\usepackage[dvips]{graphics, graphicx}

\usepackage[english]{babel}

\def\eps{\varepsilon}

\newcounter{num}[section]

\newcommand{\Th}{\refstepcounter{num}
{\bf Theorem \arabic{section}.\arabic{num} }}

\newcommand{\Lemma}{\refstepcounter{num}
{\bf Lemma \arabic{section}.\arabic{num} }}

\newcommand{\Cor}{\refstepcounter{num}
{\bf Corollary \arabic{section}.\arabic{num} }}

\newcommand{\Proof}{{\bf Proof. }}

\def\eps{\varepsilon}
\def\_phi{\varphi}

\def\d{\delta}

\def\la{\lambda}

\def\F{\widehat}

\def\f{{\mathbb F}}
\def\ov{\overline}

\def\C{{\mathbb C}}

\def\E{{\mathbb E}}

\def\Z_N{{\mathbb Z}_N}
\def\Z{{\mathbb Z}}

\def\f{{\mathbb F}}
\def\Gr{{\mathbf G}}

\def\E{\mathsf E}

\author{\sc Tomasz Schoen\footnote{The author is partially supported by MNSW grant N N201
543538.} ~ and ~ Ilya~D.~Shkredov\footnote{The author is supported
Pierre Deligne's grant based on his 2004 Balzan prize, President's
of Russian Federation grant N МК--1959.2009.1, grant RFFI N
06-01-00383 and grant Leading Scientific Schools No. 691.2008.1}}
\title{On a question of Cochrane and Pinner concerning  multiplicative
subgroups}

\date{}
\begin{document}
\maketitle

\begin{center}
 Abstract
\end{center}

{\it \small
    Answering a question of  Cochrane and  Pinner, we prove that
    for any $\eps>0$, sufficiently large prime number $p$ and an arbitrary
    multiplicative subgroup $R$ of the field $\f_p^*$, $p^{\eps} \le |R| \le p^{2/3-\eps}$ the following
    holds $|R \pm R| \ge |R|^{\frac{3}{2}+\d}$, where $\d>0$ depends on $\eps$ only.
}
\\

\refstepcounter{section}
\label{sec:introduction}

{\bf \arabic{section}. Introduction.}

\bigskip

Let $p$ be a prime number, and $R\subseteq \f_p^*$ be a
multiplicative subgroup. Such subgroups were studied by various
authors,  see e.g. \cite{Bou_prod1}--\cite{K_Tula},
\cite{Shkredov_RplusR}, \cite{Yekhanin_subgroups}.
Heath--Brown and  Konyagin  \cite{Heath_B-K} proved that for any
multiplicative subgroup $R\subseteq \f_p^*$ with  $|R|=O(p^{2/3})$
we have  $|R\pm R| \gg |R|^{3/2}$. Obviously,  the result is best
possible for subgroups of size approximately $p^{2/3}$. On the other
hand  Cochrane and  Pinner asked about the possibility of improving
the last bound for smaller subgroups (see \cite{Waring_Z_p},
Question 2). The aim of this paper is to  answer their question in
the affirmative. Our main result can be stated as follows.

\bigskip

\Th
{\it
    Let $\eps \in (0,1]$ be a real number.
    Then there exists a positive integer $p_0 (\eps)$ and a real $\d \in (0,1]$, $\d=\d(\eps)$
    such that for all primes $p \ge p_0$
    and every multiplicative subgroup $R\subseteq \f_p^*$, $p^{\eps} \le |R| \le c p^{2/3}$
    the following holds
    \begin{equation}\label{f:main_7R_1}
        |R \pm R| \gg \min\{ |R| p^{\frac13} (\log |R|)^{-\frac13},  |R|^{\frac{3}{2}} p^{\delta} \} \,.
    \end{equation}
    Furthermore, for every $|R| \le cp^{2/3}$, we have
    \begin{equation}\label{f:main_7R_1+}
        |R \pm R| \gg \min \big\{ |R| p^{\frac13} (\log |R|)^{-\frac13},
                            \max\{  |R|^{\frac{7}{3}} p^{-\frac{1}{3}} (\log |R|)^{-\frac23},
                                        |R|^{\frac{27}{14}} p^{-\frac{1}{7}} (\log |R|)^{-\frac47} \} \big\} \,.
    \end{equation}
}
\label{t:main_7R_1}

\bigskip

The formula (\ref{f:main_7R_1}) can be applied for small subgroups
and the inequality (\ref{f:main_7R_1+}) gives good bounds for large
subgroups. In particular $|R|^{{27}/{14}} p^{-{1}/{7}} (\log
|R|)^{-1/2} \gg |R|^{3/2}$, provided that $|R|>p^{1/3+\eps}$,
$\eps>0$ and if $|R|\sim p^{1/2}$ then (\ref{f:main_7R_1+}) yields
$|R\pm R|\gg |R|^{5/3-\eps}.$

In \cite{Glibichuk_zam}  Glibichuk proved the following interesting
result (see also \cite{R}).

\bigskip

\Th
{\it
    Let $A,B \subseteq \f_p$ be two sets, $|A| |B| > p$.
    Suppose that $B=-B$ or $B\cap (-B) = \emptyset$.
    Then $8AB = \f_p$.
}
\label{t:8AB}

\bigskip


\bigskip

\Cor
{\it
    Let $R\subseteq \f_p^*$ be a multiplicative subgroup and $|R| > \sqrt{p}$.
    Then $8R = \f_p$.
}
\label{cor:8R}

\bigskip



We derive from Theorem \ref{t:main_7R_1} that for all sufficiently
large $p$, in Corollary \ref{cor:8R} one can replace $8R$ by $6R$
(for precise formulation see Theorem \ref{t:main_7R_new} below). The
last result is connected with
an intriguing  question concerning  basis properties of
multiplicative subgroups. Let $R \subseteq \f_p^*$ be such a
subgroup, $|R| \ge p^{1/2+\eps}$, $\eps>0$. What is the least $l$
such that $lR$ contains $\f_p^*$? A well--known  hypothesis states
that $l$ equals two. As was showed in \cite{Yekhanin_subgroups} the
last hypothesis is not true in general finite fields ${\f}_{p^n}$,
$n\to \infty$ even for subgroups $R$ with restriction $|R| \le
|\f_{p^n}|^{2/3-\eps}$, $\eps>0$.

The main idea of our proof was used in various recent papers
\cite{Katz_Koester},
\cite{Sanders_non-abealian}---\cite{Schoen_Freiman} (see formula
(\ref{f:intr_main_f}) below) and can be described as follows. Let
$A$ be a subset of an abelian group $\Gr$. A common argument, which
can be found in many proofs of additive results uses the estimate of
energy  $\E(A)\ge |A|^4/|S|,$ where $S=A-A.$    We show, roughly
speaking, that if $\E(A)$ does not exceed  $|A|^{4-\eps}/|S|$ (but
in a stronger sense, actually we use higher moments of convolution,
see Corollary \ref{cor:huge_A+A_s}) then typically $|A - A_s|\gg
|S|^{1-c\eps},$ where $A_s = A \cap (A-s)$.
Then, using an obvious  inequality
\begin{equation}\label{f:intr_main_f}
    |S \cap (S-s)| \ge |A - A_s| \,
\end{equation}
we prove  that $\E(S)\gg |S|^{3-c'\eps}.$

Now, assume that  $A$ is a multiplicative subgroup satisfying $|S|
\ll |A|^{3/2+\eps}$ (the case $|A+A| \ll |A|^{3/2+\eps}$ is very
similar). From a result of Heath--Brown and Konyagin
\cite{Heath_B-K} it follows immediately that   $\E(A)\ll
|A|^{4+\eps}/|S|.$ Therefore, we have
$\E(S)\gg |S|^{3-c'\eps}.$
However, the last inequality cannot be true provided that $\eps$ is
small enough. It is easy to observe that
the set
$S$ just a union of some cosets and we know that
each of the coset is
uniformly distributed
(see Theorem \ref{t:BGK} and Corollary \ref{cor:pre_Fourier_R} below).
The obtained contradiction
completes
the proof.

We conclude with few comments regarding the notation used in this
paper.
All logarithms used in the paper are to base $2.$ By  $\ll$ and
$\gg$ we denote the usual Vinogradov's symbols. Finally, with a
slight abuse of notation we use the same letter to denote a set
$S\subseteq \Gr$ and its characteristic function $S:\Gr\rightarrow
\{0,1\}.$

The second author is grateful to S.V. Konyagin for useful
discussions.

\bigskip

\refstepcounter{section}
\label{sec:previous}

{\bf \arabic{section}. Previous results.}
\bigskip

In this  section we collect basic definitions we shall use later on
and quickly recall known  additive properties of multiplicative
subgroups.

Let $\Gr$ be a finite Abelian group, $N=|\Gr|.$
It is well--known~\cite{Rudin_book} that the dual group $\F{\Gr}$ is isomorphic to $\Gr.$
Let $f$ be a function from $\Gr$ to $\C.$  We denote the Fourier transform of $f$ by~$\F{f},$
\begin{equation}\label{F:Fourier}
  \F{f}(\xi) =  \sum_{x \in \Gr} f(x) e( -\xi \cdot x) \,,
\end{equation}
where $e(x) = e^{2\pi i x}$.
We rely on the following basic identities
\begin{equation}\label{F_Par}
    \sum_{x\in \Gr} |f(x)|^2
        =
            \frac{1}{N} \sum_{\xi \in \F{\Gr}} \big|\widehat{f} (\xi)\big|^2 \,.
\end{equation}
\begin{equation}\label{svertka}
    \sum_{y\in \Gr} \Big|\sum_{x\in \Gr} f(x) g(y-x) \Big|^2
        = \frac{1}{N} \sum_{\xi \in \F{\Gr}} \big|\widehat{f} (\xi)\big|^2 \big|\widehat{g} (\xi)\big|^2 \,.
\end{equation}
If
$$
    (f*g) (x) := \sum_{y\in \Gr} f(y) g(x-y) \quad \mbox{ and } \quad (f\circ g) (x) := \sum_{y\in \Gr} f(y) g(y+x)
$$
 then
\begin{equation}\label{f:F_svertka}
    \F{f*g} = \F{f} \F{g} \quad \mbox{ and } \quad \F{f \circ g} = \ov{\F{\ov{f}}} \F{g} \,.
\end{equation}
Let also $f^c (x):= f(-x)$ for any function $f:\Gr \to \C$.

Write $\mathsf E(A,B)$ for {\it additive energy} of two sets $A,B
\subseteq \Gr$ (see e.g. \cite{Tao_Vu_book}), that is
$$
    \E(A,B) = |\{ a_1+b_1 = a_2+b_2 ~:~ a_1,a_2 \in A,\, b_1,b_2 \in B \}| \,.
$$
If $A=B$ we simply write $\E(A)$ instead of $\E(A,A).$ Clearly,
\begin{equation}\label{f:energy_convolution}
    \E(A,B) = \sum_x (A*B) (x)^2 = \sum_x (A \circ B) (x)^2 = \sum_x (A \circ A) (x) (B \circ B) (x)
    \,,
\end{equation}
and by (\ref{svertka}),
\begin{equation}\label{f:energy_Fourier}
    \E(A,B) = \frac{1}{N} \sum_{\xi} |\F{A} (\xi)|^2 |\F{B} (\xi)|^2 \,.
\end{equation}

Now let $\Gr=\f_p$, where $p$ is a prime number. We call a set
$Q\subseteq \f_p\, R$--invariant if  $QR=Q.$ First of all let us
estimate Fourier coefficients of an arbitrary $R$--invariant set.

\bigskip

\Lemma
{\it
    Let $R \subseteq \f_p^*$ be a multiplicative subgroup and let $Q$ be a nonempty $R$--invariant set.
    Then for all $\xi \neq 0$ the following holds
    \begin{equation}\label{}
        |\F{Q} (\xi) | <  |Q|^{1/2}p^{1/2}{|R|}^{-1/2} \,.
    \end{equation}
}
\label{l:R-invariant}
\Proof
By Parseval identity and $R$--invariance
$$
    |R| |\F{Q} (\xi)|^2 \le \sum_{\xi \neq 0} |\F{Q} (\xi)|^2 = p |Q| - |Q|^2 < p |Q|
$$
and the result follows. $\hfill\Box$

\bigskip

Using Stepanov's method,  Heath-Brown and  Konyagin proved the
following theorem (see \cite{Heath_B-K}).

\bigskip

\Th
{\it
    Let $R \subseteq \f_p^*$ be a multiplicative subgroup and let  $Q \subseteq \f_p^*$ be a $R$--invariant set such that
    that $|Q| \ll \frac{p^3}{|R|^3}$.
    Then
    \begin{equation}\label{f:Stepanov_2/3}
        \sum_{\xi \in Q} (R \circ R) (\xi)
            \ll
                |R| |Q|^{2/3} \,.
    \end{equation}
}
\label{t:Stepanov_2/3}

\bigskip

We shall apply  a lemma from \cite{Heath_B-K} (see also
\cite{KS1,K_Tula}) which is
a
 consequence of the theorem above.

\bigskip

\Lemma
{\it
    Let $R \subseteq \f_p^*$ be a multiplicative subgroup, $|R| = O(p^{2/3})$.
    Then $$ \E(R) \ll |R|^{5/2} \,.$$
}
\label{l:t^5/2}

\bigskip

By Cauchy--Schwarz inequality it follows immediately  that $|R\pm R|
\gg |R|^{3/2}$. We shall use another consequence of Lemma
\ref{l:t^5/2}.

\bigskip

\Lemma
{\it
    Let $R \subseteq \f_p^*$ be a multiplicative subgroup, $|R| = O(p^{2/3})$,
    and let $Q$ be a nonempty $R$--invariant set.
    Then for all $\xi \neq 0$ the following holds
    \begin{equation}\label{}
        |\F{Q} (\xi) | \ll |Q|^{3/4} p^{1/4} |R|^{-3/8} \,.
    \end{equation}
} \label{l:R-invariant_new} \Proof We use the fact that $Q$ is a
disjoint union of some cosets $x_j R$, $j = 1,2,\dots, |Q|/|R|$. By
H\"{o}lder inequality, $R$--invariance and Lemma \ref{l:t^5/2}, we
have
\begin{eqnarray*}
    |\F{Q} (\xi)|
        &\le&
            \sum_{i=1}^{|Q|/|R|} |\F{R} (x_i \xi)|
                \le
                    \Big( \sum_{i=1}^{|Q|/|R|} |\F{R} (x_i \xi)|^4 \Big)^{1/4}
                    (|Q|/|R|)^{3/4}\\
                    &=&  |Q|^{3/4} |R|^{-1} \Big( \sum_{i=1}^{|Q|/|R|} |R| |\F{R} (x_i \xi)|^4 \Big)^{1/4}
                                =|Q|^{3/4} |R|^{-1} \Big( \sum_{x} |\F{R} (x)|^4
                                \Big)^{1/4}\\
                    &\le& |Q|^{3/4} |R|^{-1} (p\, \E(R) )^{1/4} \ll
                                                |Q|^{3/4} p^{1/4} |R|^{-3/8}
\end{eqnarray*}
and the result follows. $\hfill\Box$

\bigskip

Also Lemma \ref{l:R-invariant}, Lemma \ref{l:R-invariant_new} and
Lemma \ref{l:t^5/2} implies the following upper estimate for Fourier
coefficients of a multiplicative subgroup (see \cite{KS1} or
\cite{K_Tula}). For completeness we recall the proof.

\bigskip

\Cor
{\it
    Let $R \subseteq \f_p^*$ be a multiplicative subgroup.
    Then $\max_{\xi \neq 0} |\F{R} (\xi)| < \sqrt{p}$.
    Suppose, in addition, that $|R| = O(p^{2/3})$.
    Then
    \begin{equation}\label{f:previous_F_bound}
        \max_{\xi \neq 0} |\F{R} (\xi)| \ll \min\{ p^{\frac{1}{4}} |R|^{\frac{3}{8}}, p^{\frac{1}{8}} |R|^{\frac{5}{8}} \} \,.
    \end{equation}

} \label{cor:pre_Fourier_R} \Proof Let $\rho = \max_{x \neq 0}
|\F{R} (x)|$. The first estimate $\rho < \sqrt{p}$\, is a
consequence on Lemma \ref{l:R-invariant}. We should check
(\ref{f:previous_F_bound}). To obtain the first bound in the
formula, we apply Lemma \ref{l:R-invariant_new} with $Q=R$. Further,
let $\xi \neq 0$ be an arbitrary residual. By $R$--invariance and
H\"{o}lder inequality we have
\begin{equation*}\label{f:pre_4_moment}
       |R|^2 |\F{R} (\xi)|^2 = \Big|\sum_{x,y \in R} e( -\xi xy)\Big|^2
            \le
                |R| \sum_{x\in R} |\F{R} (\xi x)|^2 = |R| \sum_{x} (\xi R * (\xi R)^c) (x) \F{R^c}(x) \,.
\end{equation*}
Thus
\begin{eqnarray*}\label{f:pre_4_moment'}
    |R|^4 |\F{R} (\xi)|^4
                    &\le& |R|^2 \sum_{x} (\xi R*(\xi R)^c) (x) \cdot \sum_{x} (\xi R* (\xi R)^c) (x) |\F{R}
                    (x)|^2 \nonumber \\
                    &=& |R|^4 \sum_{\xi} (\xi R*(\xi R)^c) (\xi) |\F{R} (x)|^2 \,,
\end{eqnarray*}
and
\begin{equation*}\label{f:pre_4_moment''}
    |R|^8 |\F{R} (\xi)|^8
                    \le
                        |R|^8 \E(R)^2 p \,.
\end{equation*}
Now the assertion follows from  Lemma \ref{l:t^5/2}. $\hfill\Box$

\bigskip

From the above  estimates of Fourier coefficients of subgroups one
can deduce its basis properties.

\bigskip

\Th
{\it
    Let $R\subseteq \f_p^*$ be a multiplicative subgroup, $l\ge 4$ be a positive integer.
    Suppose that $|R| \gg \min\{ p^{\frac{2l+2}{5l-3}}, p^{\frac{l+5}{3l+3}} \}$.
    Then $\f^*_p \subseteq lR$.
    If $l\ge 2$ and $|R|>p^{\frac{l+1}{2l}}$ then $\f^*_p \subseteq lR$.
}
\label{t:pre_main_7R}





\Proof  Suppose that $lR \nsupseteq \f^*_p$. Then for some $\la\neq
0$ we have
\begin{equation}\label{tmp:19.04.2010_1}
    0 = \sum_{x} \F{R}^l (x) \F{\la R} (x) = |R|^{l+1} + \sum_{x \neq 0} \F{R}^l (x) \F{\la R} (x) \,.
\end{equation}
Using Corollary \ref{cor:pre_Fourier_R} and Parseval formula to
estimate the second term in (\ref{tmp:19.04.2010_1}), we get the
required result. $\hfill\Box$

\bigskip

In particular, if $l=7$ and $|R| \gg \sqrt{p}$ then $\f^*_p
\subseteq lR$.
For large $l$ better bounds are known (see \cite{KS1,BGK}).

Finally, we recall a well--known result of Bourgain, Glibichuk and Konyagin \cite{BGK} on
Fourier coefficients of
multiplicative subgroups.

\bigskip

\Th
{\it
    Let $\eps \in (0,1]$ be a real number.
    Then there exists a positive integer $p_0 (\eps)$ and a real $\eta \in (0,1]$
    such that for all primes $p \ge p_0$
    and any multiplicative subgroup $R\subseteq \f_p^*$, $|R| \ge
    p^{\eps}$ we have
    \begin{equation}\label{f:BGK}
        \max_{\la \in \f_p^*}\, \Big|\sum_{x\in R} e(\la x) \Big| \le |R| p^{-\eta} \,.
    \end{equation}
}
\label{t:BGK}

Notice that for large subgroups (i.e. $|R|>p^{1/4+\eps}$, $\eps>0$)
better estimates hold (see Corollary \ref{cor:pre_Fourier_R} above
and \cite{Heath_B-K,KS1,K_Tula}). Finally, we remark that the
condition $|R| \ge p^{\eps}$ in the theorem above can be replaced by
a weaker inequality $|R| \gg \exp (C' \log p /\log \log p)$, where
$C'>1$ is a suitable  constant (see \cite{Bourgain_new_sum-prod}).
\bigskip

\refstepcounter{section} \label{sec:addcom}
\bigskip

{\bf \arabic{section}. Additive combinatorics.}

\bigskip

We return for a moment to a general case of an arbitrary Abelian
group $\Gr$. For any set $A \subseteq \Gr$ and any element $s\in
A-A$ define  the set $A_s = A \cap (A-s)$
(see papers \cite{Katz_Koester}---\cite{Schoen_Freiman}). Clearly,
$A_s \neq \emptyset$ and $|A_s| = (A\circ A) (-s) = (A\circ A) (s)$.
Furthermore, set
$$\E_3(A)=\sum_s (A\circ A) (s)^3$$
and
$$\E_4(A)=\sum_s (A\circ A)(s)^4.$$
The following simple lemma exhibits interesting relations between
energies of  $A_s$ and the quantities $\E_3(A), \E_4(A).$

\bigskip

\Lemma {\it
    Let $\Gr$ be an Abelian group.
    For every set $A\subseteq \Gr$ we have
    \begin{equation*}\label{f:E(A_s,A)}
        \sum_{s\in A-A} \E(A,A_s) = \E_3(A) \,
    \end{equation*}
    and
    \begin{equation*}\label{f:E(A_s,A_t)}
        \sum_{s\in A-A} \sum_{t\in A-A} \E(A_s,A_t) = \E_4(A) \,.
    \end{equation*}
} \label{l:E(A_s,A)} \Proof
   Observe that for every $x,w,s \in \Gr$
$$
    A_s (x) A_s (x+t) = A(x) A(x+s) A(x+t) A(x+s+t) = A_t (x) A_t (x+s) \,.
$$
Summing over $x$, we get
\begin{equation*}\label{tmp:1.08.2010_1}
    (A_s \circ A_s) (t) = (A_t \circ A_t) (s) \,, \quad \quad s,t\in \Gr \,.
\end{equation*}
Clearly,
$$
    \sum_{s} (A_t \circ A_t) (s) = |A_t|^2 = (A\circ A)(t)^2 \,.
$$
Thus, by (\ref{tmp:1.08.2010_1}) we have
$$
 \sum_{s} \E(A,A_s)=\sum_s\sum_t (A \circ A) (t)(A_s \circ A_s) (t)= \sum_t  (A \circ A) (t)\sum_s (A_t \circ A_t) (s)=\E_3(A)\,.
$$
Similarly,
\begin{eqnarray*}
  \sum_{s} \sum_{t} \E(A_s,A_t) &=& \sum_s\sum_t \sum_x(A_s \circ A_s) (x)(A_t \circ A_t) (x)\\
  &=& \sum_x \sum_s (A_x \circ A_x) (s)\sum_t(A_x \circ A_x) (t)=\E_4(A) \,,
\end{eqnarray*}
which proves  Lemma \ref{f:E(A_s,A_t)}. $\hfill\Box$

\bigskip

Next, we show that if $\E_3(A)$ and $\E_4(A)$ are small then
typically $|A-A_s|$ and $|A_s-A_t|$ are large, respectively.

\bigskip

\Cor {\it
    Let $A$ be a subset of an abelian group $\Gr$ with at least $2$ elements.
    Then
    \begin{equation}\label{f:cor_huge_A+A_s}
        \sum_{s\neq 0} |A\pm A_s| \ge 2^{-1} |A|^6 \E_3(A)^{-1}  \,
    \end{equation}
    and
    \begin{equation}\label{f:cor_huge_A_s+A_t}
        \sum_{s\neq 0} \sum_{t\neq 0} |A_s\pm A_t| \ge  4^{-1} |A|^8 \E_4(A)^{-1} \,.
    \end{equation}
    In particular there exist $s\neq 0$ such that  $|A-A_s| \ge 2^{-1}|A|^6\E_3(A)^{-1}|A-A|^{-1}$
     and $s,t\neq 0$ such that $|A_s-A_t| \ge  4^{-1}|A|^8\E_4(A)^{-1}|A-A|^{-2}.$
} \label{cor:huge_A+A_s}

\Proof Fix $s\in A-A$. By Cauchy--Schwarz inequality and
(\ref{f:energy_convolution}), we obtain
$$
    ( |A| |A_s| )^2 = \Big( \sum_z (A\circ A_s) (z) \Big)^2 = \Big( \sum_z (A * A_s) (z) \Big)^2
        \le \E(A,A_s) \cdot |A\pm A_s| \,,
$$
so that
\begin{equation*}\label{f:intermediate_0.5}
    |A| (|A|^2-|A|)
        = \sum_{s\in (A-A) \setminus \{ 0 \}}  |A| |A_s|
          \le  \sum_{s\in (A-A) \setminus \{ 0 \}} \E(A,A_s)^{1/2} \cdot |A \pm A_s|^{1/2} \,.
\end{equation*}
Applying once again   Cauchy--Schwarz inequality  and Lemma
\ref{l:E(A_s,A)}, we have
\begin{eqnarray*}
    2^{-1} |A|^3
        &\le&
            \sum_{s\neq 0} \E (A,A_s)^{1/2}  |A \pm A_s|^{1/2}
                \le \Big( \sum_{s\neq 0} \E (A,A_s) \Big)^{1/2}  \Big( \sum_{s\neq 0} |A \pm A_s| \Big)^{1/2}\\
         &\le& \E_3 (A)^{1/2}  \Big( \sum_{s\neq 0} |A \pm A_s| \Big)^{1/2} \,.
\end{eqnarray*}
The second inequality one can prove using very similar argument.
$\hfill\Box$

\bigskip

We  bound from the above  $\E_3(R)$ and $\E_4(R)$, where $R$ is a
multiplicative subgroup.

\bigskip

\Lemma {\it
    Let $R \subseteq \f_p^*$ be a multiplicative subgroup
    and $|R| \ll p^{2/3}$.
    Then
        \begin{equation}\label{f:E_*(R)_2}
            \E_3 (R) \ll |R|^3 \log |R|  \text{ ~  and ~  }  \E_4 (R) \le |R|^4+ O(|R|^{11/3})\,.
        \end{equation}
} \label{l:E_*(R)} \Proof
Let $a$ be a parameter, and put $n=|R|$. We have
$$
    \E_3 (R) \le a^2 n^2 + \sum_{s ~:~ (A\circ A) (s) \ge a} (A\circ A)(s)^3 \,.
$$
Let us arrange values $(A\circ A) (s)$, $s\in \f_p / R$ in
decreasing order and denote its values as $N_1 \ge N_2 \ge \dots$.
By Theorem \ref{t:Stepanov_2/3}, we have $N_j \ll n^{2/3} j^{-1/3}$.
Hence
$$
    \E_3 (R) \le a^2 n^2 + n \sum_{j ~:~ j\ll n^2 / a^3} N^3_j
        \ll
            a^2 n^2 + n^3 \sum_{j ~:~ j\ll n^2 / a^3} j^{-1}
                \ll
                    a^2 n^2 + n^3 \log \frac{n^2}{a^3} \,.
$$
Taking $a=n^{1/2}$ we obtain the required result. To establish the
second estimate, one similarly  bounds the contribution of terms
with $s\neq 0$ to $\E_4(R)$ by $\ll |R|^{11/3}.$ $\hfill\Box$

\bigskip

{\bf Remark.} Suppose that $R$ is a multiplicative subgroup, $|R|
\ll p^{2/3}$. Clearly,  by (\ref{f:energy_Fourier}) we have $\E(A)
\ge p^{-1} |A|^4$ for every set $A\subseteq \Gr$. Hence, by
Corollary \ref{cor:huge_A+A_s} and Lemma \ref{l:E_*(R)}, we can find
$s\in R-R$, $s\neq 0$ such that $|R \pm R_s| \gg \frac{|R|^3}{p\log
|R|}$. It gives interesting consequences for subgroups $R$
satisfying
$|R| \ge c p^{2/3}$, where $c>0$ is a small constant. Indeed, in
this case there is $s\in R-R$, $s\neq 0$
with
 $|R \pm R_s| \gg \frac{p}{\log |R|}$.
Actually,
from (\ref{f:cor_huge_A+A_s}) one can deduce that there are at least
$\gg \frac{|R|^3}{p \log |R|} \sim \frac{p}{\log |R|}$ such $s$.
Notice that the restriction $|R| \ge cp^{2/3}$ implies
$|R-R| \gg p$ (see e.g. \cite{Shkredov_RplusR}). Furthermore, from
Corollary \ref{cor:huge_A+A_s} and  Lemma \ref{l:E_*(R)} it follows
immediately that for some $s,t\neq 0$ we have $|R_s\pm R|\gg
p^{2/3}.$ However, in this case we can prove even more.
Let $P\subseteq (R-R)\setminus \{0\}$ be the set of all popular
differences, i.e. $s\in \f_p^*$ having at least $ |R|^2/(10p)\gg
p^{1/3}$ representations as $r_1-r_2,\, r_1,r_2\in R$. By H\"older
inequality, (\ref{f:cor_huge_A_s+A_t}) and  (\ref{f:E_*(R)_2})  we
have
\begin{eqnarray*}|R|^4&\ll &
\sum_{s,t\in P}|R_s||R_t|\le \big(\sum_{s,t\in P}\E
(R_s,R_t)\big)^{1/2}\big(\sum_{s,t\in
P}(|R_s||R_t|)^{2}\big)^{1/4}\Big(\sum_{s,t\in
P}\frac{|R_s\pm R_t|}{|R_s||R_t|}\Big )^{1/4}\\
&\le & \E_4(R)^{1/2}\E(R)^{1/2}\Big(\sum_{s,t\in P}\frac{|R_s\pm
R_t|}{|R_s||R_t|}\Big )^{1/4}\le |R|^{13/4}\Big(\sum_{s,t\in
P}\frac{|R_s\pm R_t|}{|R_s||R_t|}\Big )^{1/4},
\end{eqnarray*}
hence
$$|R|^{3}\ll \sum_{s,t\in P}
\frac{|R_s\pm R_t|}{|R_s||R_t|}.$$ Thus, if $|R|\gg p^{2/3}$ there
are $\gg p^2$ pairs $s,t$ such that $|R_s|,|R_t|\gg p^{1/3}$ and
$|R_s-R_t|\gg |R_s||R_t|.$



\refstepcounter{section} \label{sec:proof}
\bigskip

\bigskip

{\bf \arabic{section}.  Proof of Theorem 1.1.}

\bigskip

Let $S=R-R$, $S'=R+R $, $n=|R|$, and $m = |S|$.
We may assume that $n \ge 2$.
Observe  that the sets
$ S\setminus \{0\}$ and  $S'\setminus \{0\} $ are $R$--invariant.
Now we apply  arguments used
\cite{Katz_Koester},
\cite{Sanders_non-abealian}---\cite{Schoen_Freiman} Fix $s\in S$. In
view of  $(A - A_s) \subseteq S$, $(A - A_s)  \subseteq S-s$ and $(A
+ A_s) \subseteq S'$, $(A + A_s) \subseteq S'-s$, we have
\begin{equation*}\label{f:sumset_property}
    (S\circ S) (s) \ge |A- A_s|\,, \quad (S'\circ S') (s) \ge |A + A_s|\,, \quad \quad s \in S \,.
\end{equation*}
We will deal with the sign "-" and the set $S$.
The case of $S'$  can be handle in the same way.
By Corollary \ref{cor:huge_A+A_s} and  (\ref{f:sumset_property})
$$
     n^6 \E_3 (R)^{-1} \ll \sum_s S (s) (S\circ S) (s) \,.
$$
Using Fourier transform and Parseval formula, we infer that
$$
     n^6 \E_3 (R)^{-1} p \ll \sum_{\xi} |\F{S} (\xi)|^2 \F{S} (\xi) \le m^3 +  \rho m p \,,
$$
where $\rho=\max_{\xi \neq 0} |\F{S} (\xi)|$. Whence either $n^6
\E_3(R)^{-1} p \ll  m^3$ or $n^6 \E_3 (R)^{-1} \ll   \rho m$. In the
first case, by Lemma \ref{l:E_*(R)}, we have
\begin{equation*}\label{f:E_1}
   n p^{1/3} (\log n)^{-1/3}\ll  m \,.
\end{equation*}
Now suppose that $n^6 \E_3 (R)^{-1} \ll   \rho m$.
By Theorem \ref{t:BGK}  $\rho \le m p^{-\eta}$. Therefore, by Lemma
\ref{l:E_*(R)}, we get
\begin{equation*}\label{f:E_2}
    {n^3}({\log n})^{-1} \ll m^2 p^{-\eta}\, ,
\end{equation*}
which proves (\ref{f:main_7R_1}). Here one can take any $\delta <
\eta/2$.

Next,  we  prove (\ref{f:main_7R_1+}). We can consider  only the
case $n^6 \E_3 (R)^{-1} \ll  \rho m$. Using Lemma
\ref{l:R-invariant}, we have $\rho < ( p m/n)^{1/2}+1$. Hence
$$
    {n^3}({\log n})^{-1} \ll m^{3/2} p^{1/2}{n}^{-1/2}\,,
$$
so that $m \gg n^{{7}/{3}} p^{-{1}/{3}} (\log n)^{-2/3}$. Finally,
applying Lemma \ref{l:R-invariant_new}, we obtain
$$
    {n^3}({\log n})^{-1} \ll  m^{7/4} p^{1/4} n^{-3/8}\,,
$$
which gives
 $m \gg n^{{27}/{14}} p^{-{1}/{7}} (\log
n)^{-4/7}$. This completes the proof. $\hfill\Box$


\bigskip

Our main result allows us to  refine Theorem \ref{t:pre_main_7R}.

\bigskip

\Th
{\it
    Let $R\subseteq \f_p^*$ be a multiplicative subgroup and $\kappa > \frac{41}{83}.$
    Suppose that $|R| \ge p^{\kappa}$,
    then for all sufficiently large $p$ we have $\f^*_p \subseteq 6R$.
}
\label{t:main_7R_new}
\\
\Proof Suppose that $R\subseteq \f_p^*$ is a multiplicative subgroup
with  $n=|R| \ge p^{\kappa}$, $\kappa > \frac{41}{83}$. Furthermore,
put $S = R+R$, $n=|R|$, $m = |S|$, and $\rho = \max_{\xi \neq 0}
|\F{R} (\xi)|$.  If $6R\not\subseteq \f_p^*$ then for some
$\lambda\neq 0,$ we have
$$ 0 = \sum_{\xi} \F{S}^2 (\xi) \F{R}^2 (\xi) \F{\la R} (\xi)=m^2n^3+\sum_{\xi\neq 0}
\F{S}^2 (\xi) \F{R}^2 (\xi) \F{\la R} (\xi).$$ Therefore, by   the
second inequality of (\ref{f:previous_F_bound})
and Parseval identity we get
$$n^3m^2\le \rho^3mp\ll (p^{1/8}n^{5/8})^3mp.$$
Now applying the second inequality of (\ref{f:main_7R_1+}),  $m \gg
n^{{7}/{3}} p^{-{1}/{3}} (\log n)^{-2/3},$ we obtain the required
result. $\hfill\Box$

\bigskip

Finally, let us remark that using   similar argument as  in the
proof of Theorem \ref{t:main_7R_new} one can show that $|4R|\ge
p/2,$ provided that $|R|\ge p^{\kappa},$ where $\kappa> 11/23,$ and
$p$ is large enough.




\end{document}